\documentstyle[12pt,draft]{article}

 \newcounter{abceqn}

 \newcounter{abcfig}

\newcommand{\be}{\beta}

\newcommand{\Om}{\Omega}

\newcommand{\pa}{\partial}

\newcommand{\ga}{\gamma}

\newcommand{\dl}{\delta}
\newcommand{\Dl}{\Delta}

\newcommand{\al}{\alpha}

\newcommand{\la}{\lambda}

\renewcommand{\thesection}{\arabic{section}.}

\newcommand{\eqnsection}[1]{
	\section{#1}
	\setcounter{equation}{0}
	\renewcommand{\theequation}{\thesection\arabic{equation}}
	\setcounter{figure}{0}
	\renewcommand{\thefigure}{\arabic{figure}}
	\setcounter{remark}{0}
	\renewcommand{\theremark}{\thesection\arabic{remark}}
	\setcounter{theorem}{0}
	\renewcommand{\thetheorem}{\thesection\arabic{theorem}}
	\setcounter{lemma}{0}
	\renewcommand{\thelemma}{\thesection\arabic{lemma}}
}

\title{{\bf A Lax Pair for the 2D Euler Equation}}

\author{ \\ \\ \\ \\ 
Yanguang (Charles)\ \ Li  \thanks{This work is supported 
by the Guggenheim Fellowship.}
\\  \\  \\ Department of Mathematics, 
 \\ \\ University of Missouri \\ \\ 
Columbia, MO 65211 \\ \\ E-mail: cli@math.missouri.edu \\ }

\date{\today}

\begin{document}
\bibliographystyle{unsrt}
\maketitle
\newpage
\begin{abstract}    
A Lax pair for the 2D Euler equation is found. \\



\end{abstract}

\newtheorem{lemma}{Lemma}
\newtheorem{theorem}{Theorem}
\newtheorem{corollary}{Corollary}
\newtheorem{remark}{Remark}
\newtheorem{definition}{Definition}
\newtheorem{proposition}{Proposition}
\newtheorem{assumption}{Assumption}

\eqnsection{A Lax Pair for the 2D Euler Equation}

This is to report that a Lax pair for the 2D Euler equation is found.
We write the 2D Euler equation in the vorticity form,
\begin{equation}
{\pa \Om \over \pa t} + \{ \Psi, \Om \} = 0 \ ,
\label{euler}
\end{equation}
where $\Om$ is the vorticity, $\Psi$ is the stream function, and 
the bracket $\{\ \}$ is defined as
\[
\{ f, g\} = (\pa_x f) (\pa_y g) - (\pa_y f) (\pa_x g) \ .
\]
Let us denote the $x$-directional and the $y$-directional velocities
by $u$ and $v$ respectively. Then
\[
u = - {\pa \Psi \over \pa y}\ , \ \ v = {\pa \Psi \over \pa x}\ ,
\ \ \Om = {\pa v \over \pa x} - {\pa u \over \pa y}\ , \ \ 
\Dl \Psi = \Om \ .
\]
The Lax pair is given as
\begin{equation}
\left \{ \begin{array}{l} 
L \varphi = \la \varphi \ ,
\\
\pa_t \varphi + A \varphi = 0 \ ,
\end{array} \right.
\label{laxpair}
\end{equation}
where
\[
L \varphi = \{ \Om, \varphi \}\ , \ \ \ A \varphi = \{ \Psi, \varphi \}\ ,
\]
and $\la$ is a complex constant, and $\varphi$ is a complex-valued function.
The compatibility condition of the Lax pair (\ref{laxpair}) gives the 
2D Euler equation (\ref{euler}), i.e. 
\[
\pa_t L = [L, A]\ ,
\]
where $[L, A] = LA -AL$, gives the Lax representation of the 2D Euler 
equation (\ref{euler}). 

\begin{remark}
With the recent development on chaos in partial differential equations 
\cite{LM94} \cite{LMSW96} \cite{Li99}, I am interested in building a 
dynamical system theory for 2D Euler equation under periodic boundary 
condition \cite{Li00b} \cite{Li00c}. In particular, I am investigating 
the existence v.s. 
nonexistence of homoclinic structure. For such studies, it will be 
fundamentally important to find a Lax pair (if it exists) for the 2D 
Euler equation. Then I started with Vladimir Zakharov's paper \cite{Zak89}.
Zakharov proposed the Lax pair
\[
\left \{ \begin{array}{l} 
\la D_1 \varphi + \{ \Om, \varphi \} = 0  \ ,
\\
\pa_t \varphi + \la D_2 \varphi + \{ S, \varphi \}= 0 \ ,
\end{array} \right.
\]
where
\[
D_1 = \al {\pa \over \pa x} + \be {\pa \over \pa y} \ , \ \ \ 
D_2 = \ga {\pa \over \pa x} + \dl {\pa \over \pa y} \ ,
\] 
$\al , \be , \ga , \ \mbox{and} \ \dl $ are real constants, $\la$ is a complex 
constant, $S$ is a real-valued function, and $\varphi$ is a complex-valued 
function. The compatibility condition of this Lax pair gives the 
following equation instead of the 2D Euler equation,
\[
\left \{ \begin{array}{l} 
{\pa \Om \over \pa t} + \{ S, \Om \} = 0 \ ,
\\ \\
D_1 S = D_2 \Om\ .
\end{array} \right.
\]
(Notice the misprints in the English translation of the article \cite{Zak89}.)
\end{remark}
\begin{remark}
The author is also aware of the Lax pair in the inverse Cauchy-Green tensor 
variable of the Lagrangian formulations of both 2D and 3D Euler equations 
found by Susan Friedlander and Misha Vishik \cite{FV90} \cite{VF93}.
\end{remark}

\bibliography{lax}

\end{document}